\documentclass[10pt]{article}
\usepackage{amsmath}
\usepackage{amsfonts}
\usepackage{amssymb}
\usepackage{comment}
\long\def\proof#1{\removelastskip\vskip\baselineskip\relax\noindent{\it
Proof\if!#1!\else\ \ignorespaces#1\fi.\ }\ignorespaces}

\newcommand{\lgs}[2]{\mbox{$\left(\frac{#1}{#2}\right)$}}
\newcommand{\leg}[2]{\mbox{$\left(\dfrac{#1}{#2}\right)$}}

\newcommand{\Q}{{\mathbb Q}}
\newcommand{\Z}{{\mathbb Z}}

\newcommand{\z}{\zeta}

\newcommand{\Om}{\Omega}

\newcommand{\G}{\Gamma}

\DeclareMathOperator{\sign}{sign}
\DeclareMathOperator{\asin}{asin}

\newcommand{\Proof}{{\it Proof. \/}}
\newcommand{\squareforqed}{\hbox{\rlap{$\sqcap$}$\sqcup$}}
\newcommand{\qed}{\ifmmode\squareforqed\else{\unskip\nobreak\hfil
\penalty50\hskip1em\null\nobreak\hfil\squareforqed
\parfillskip=0pt\finalhyphendemerits=0\endgraf}\fi}

\newcommand{\fp}{\qed\removelastskip\vskip\baselineskip\relax}

\newtheorem{theorem}{Theorem}[section]

\newtheorem{proposition}[theorem]{Proposition}

\newtheorem{definition}[theorem]{Definition}

\newcommand{\litem}{\par\noindent\dimen0=\parindent%
    \advance\dimen0 by-4pt
               \hangindent=\dimen0\ltextindent}

\newcommand{\ltextindent}[1]{\hbox to \hangindent{#1\hss}\ignorespaces}
\newcommand{\ltextjndent}[1]{\hbox to \hangindent{#1\hss}\ignorespaces\kern-1ex}

\renewcommand{\pmod}[1]{\allowbreak\ ({\rm{mod}}\,\,#1)}

\begin{document}
\pagestyle{plain}

\title{Rational Hypergeometric Ramanujan Identities for $1/\pi^c$: Survey and Generalizations}
\author{Henri Cohen and Jes\'us Guillera}

\maketitle
\begin{abstract}
  We give a simple unified proof for \emph{all} existing rational
  hypergeometric ramanujan identities for $1/\pi$, and give a complete survey
  (without proof) of several generalizations: rational hypergeometric
  identities for $1/\pi^c$, Taylor expansions, upside-down formulas, and
  supercongruences.
\end{abstract}

\smallskip

\section{Introduction}

In a famous paper \cite{Ram}, S.~Ramanujan gave $17$ formulas for $1/\pi$.
These formulas were proved and generalized much later by numerous authors.
In the present paper, which is mainly a survey and does not claim originality,
we have several goals. First, we want to show that the $36$ rational
hypergeometric formulas for $1/\pi$ follow by specialization from a
\emph{single} general formula. Second, we give a list of all
known rational hypergeometric formulas for $1/\pi^c$ with $c\ge2$, many
unproved. Third, we will give Taylor expansions of which the $1/\pi^c$
formulas are only the constant term. Fourth, we give a list of what can be
called \emph{upside-down} formulas. Finally, we give the
\emph{supercongruences} corresponding to all the $1/\pi^c$ formulas.
A large part of the formulas of this paper apart from the initial $1/\pi$
formulas are due to the second author.

This survey is meant to be exhaustive, which means that we would appreciate
feedback from readers who are aware of formulas that are not in our list
(and evidently of errors). Note that we do not list formulas involving
algebraic (as opposed to rational) parameters, nor do we list the second
author's two-sided formulas where the sums are over $n\in\Z$ instead of
$n\ge0$, or other generalizations.

\medskip

{\bf Acknowledgment:} We heartily thank Wadim Zudilin for enlightening
conversations.

\medskip

We recall that the \emph{Pochhammer symbol} $(x)_n$ is defined by
$(x)_n=x(x+1)\cdots(x+n-1)=\G(x+n)/\G(x)$, this last formula allowing us
to define it also when $n$ is not an integer.

\begin{definition}\begin{enumerate}
\item Let $d\ge2$ be an integer. For all $n\ge0$ we define
$$R_n(d)=\prod_{\substack{1\le i\le d\\\gcd(i,d)=1}}\dfrac{(i/d)_n}{n!}\;.$$
\item A \emph{rational hypergeometric product} $H$ is a sequence of the form
$H_n=\prod_{d\in I}R_n(d)^{v_d}$ for some finite index set $I$ and
\emph{positive} exponents $v_d$. We define the \emph{degree} of $H$ by
$\deg(H)=\sum_{d\in I}\phi(d)$.
\item A \emph{rational hypergeometric Ramanujan series} is a series of the form
$$S(H,a,P)=\sum_{n\ge0}P(n)\dfrac{H_n}{a^n}\;,$$
with $P\in\Z[X]$ and $a\in\Q^*$.
\end{enumerate}\end{definition}

Some comments are in order:
\begin{enumerate}
\item We only consider coefficients which are hypergeometric \emph{products}
(as opposed to quotients). We could allow $v_d<0$, or more general
coefficients, and indeed there is a vast amount of formulas for such general
series, but we will restrict to those, although we will later give
``upside-down'' formulas where \emph{all} the $v_d$ are negative. In
the tables that we give below, we will abbreviate $\prod_{d\in I}R_n(d)^{v_d}$
as $\prod d^{v_d}$.
\item In addition, we restrict to such
products which are \emph{rational} in the hypergeometric motive sense,
in other words such that if some irreducible fraction $i/d$ occurs,
then all irreducible fractions in $]0,1[$ with denominator $d$ occur.
\item We only consider $P\in\Z[X]$ (or $P\in\Q[X]$, which is the same up
to a multiplicative constant). Ramanujan himself gave formulas where
$P$ has coefficients in a quadratic extension of $\Q$, but we will not
consider those. Similarly, we restrict to $a\in\Q^*$.
\item It is immediate to see that $H_n\sim C/n^{\deg(H)/2}$ for some constant
$C$, so the series $S$ converges for $|a|>1$, also for $a=-1$ if
$\deg(P)<\deg(H)/2$, and for $a=1$ if $\deg(P)<\deg(H)/2-1$.
In fact, in all the examples that we will see we have $\deg(P)=(\deg(H)-1)/2$
(and in particular $\deg(H)$ is odd).
\item The hypergeometric function $\sum_{n\ge0}H_nz^n$ satisfies a linear
differential equation of degree $\deg(H)$, in other words there exists
a polynomial $Q(n,z)$ of degree $\deg(H)$ in $n$ such that
$\sum_{n\ge0}Q(n,z)H_nz^n=0$. We may thus restrict to polynomials $P$ such
that $\deg(P)<\deg(H)$, since formulas with $P$ differing by a multiple of
$Q$ are trivially equivalent. We will give examples of this below.
\end{enumerate}

In fact, the only hypergeometric products that we will consider in this paper
are products of the following $R_n(d)$:
\begin{align*}
R_n(p)&=p^{-pn}\dfrac{(pn)!}{n!^p}\text{\quad for $p=2,\ 3,\ 5$\;,}\\
R_n(4)&=2^{-6n}\dfrac{(4n)!}{(2n)!n!^2}\;,\text{\quad and\quad} R_n(6)=2^{-4n}3^{-3n}\dfrac{(6n)!}{(3n)!(2n)!n!}\;.\end{align*}

\begin{definition} A (convergent) hypergeometric Ramanujan series will be
called a $1/\pi^c$-formula if its sum $S(H,a,P)$ is equal to an algebraic
number divided by $\pi^c$.\end{definition}

A search through the (abundant) literature (together with additional
personal investigations) shows that the only algebraic
numbers that occur are of the form $\sqrt{k}$ for some $k\in\Q^*$, and
we have been able to find exactly $36$ series whose sum is of the form
$\sqrt{k}/\pi$, $10$ whose sum is of the form $\sqrt{k}/\pi^2$, plus a
single example for $1/\pi^3$ due to B.~Gourevitch and two examples for
$1/\pi^4$, one due to J.~Cullen the other to Y. Zhao. In addition, there
are a number of ``divergent'' hypergeometric series for $1/\pi^c$ which we
will mention.

\section{Rational Hypergeometric Formulas for $1/\pi$}

There are (at least) three methods for proving such formulas. The first,
due to Ramanujan, is the use of elliptic functions and generalizations,
the second is the use of modular functions, and the third is to use
WZ-type summation methods. In the present paper we only use modular functions.
In fact, we will show that \emph{all} the known rational formulas for $1/\pi$
follow from a single general formula giving an identity between
complex \emph{functions}, which can then be specialized to any CM point that
we want, and in particular to CM points giving rational formulas.

\subsection{A General Identity}

An important theorem which can be found for instance in \cite{Zag} states that
any modular form (or function) $F$ of weight $k$ (say on a congruence subgroup
of $\G$), expressed locally in terms of a modular \emph{function} $h$
(of weight $0$), is a solution of a \emph{linear} differential equation of
order $k+1$ with algebraic coefficients, which can be explicitly constructed;
if in addition $h$ is a \emph{Hauptmodul}, the coefficients can be chosen
to be polynomials. Equivalently, there exists
a sequence $u(n)$ satisfying a polynomial recurrence relation such that for
$\Im(\tau)$ sufficiently large we have $F(\tau)=\sum_{n\ge0}u(n)h(\tau)^n$.
We then have the following easy result, where from now on we denote by $D$
the differential operator $D=(1/2\pi i)d/d\tau=qd/dq$:

\begin{proposition}\label{prop:main} As above, let $F$ be a modular function
  of weight $k$ on some congruence subgroup of $\G$, $h$ a modular function
  (of weight $0$), write $F(\tau)=\sum_{n\ge0}u(n)h(\tau)^n$ for $\Im(\tau)$
  sufficiently large, and finally set
  $$G^*(\tau)=\dfrac{D(F)/F}{D(h)/h}(\tau)-\dfrac{k}{4\pi\Im(\tau)D(h)/h}\;,$$
  which is nonholomorphic but modular of weight $0$. We have the following
  general formula, which can be considered as a formula for $1/\pi$:
  $$\sum_{n\ge0}(n-G^*(\tau))u(n)h(\tau)^n=\dfrac{k}{4\pi\Im(\tau)}\dfrac{F}{D(h)/h}(\tau)\;.$$
\end{proposition}

\Proof First apply $D$ to the formula expressing $F$ in terms of $h$, so that
$$\dfrac{D(F(\tau))}{D(h(\tau))/h(\tau)}=\sum_{n\ge0}nu(n)h(\tau)^n\;.$$
Thus, if we set $G=(D(F)/F)/(D(h)/h)$, the left hand side is
$GF=G\sum_{n\ge0}u(n)h^n$, so we have the identity
$$\sum_{n\ge0}(n-G(\tau))u(n)h(\tau)^n=0\;.$$
Now $D(h)/h$ is a modular function of weight $2$, but $D(F)/F$ is only
quasi-modular: $D^*(F)/F=D(F)/F-(k/(4\pi\Im(\tau)))$ is truly modular
nonholomorphic of weight $2$, Thus, we set $G^*=(D^*(F)/F)/(D(h)/h)$,
and this gives both the formula for $G^*$ and the desired identity.\fp

Now CM theory tells us that if $\tau$ is a CM point and $h(\tau)$ and
$F(\tau)$ have algebraic Fourier coefficients, then both $h(\tau)$ and
$G^*(\tau)$ will be algebraic numbers, and $(F/(D(h)/h))(\tau)$, which
has weight $k-2$, will be an algebraic number times
$\Om_{\tau}^{k-2}$, where $\Om_{\tau}$ is a suitable period, for instance
$\Om_{\tau}=\eta(\tau)^2$. Thus, it will itself be algebraic
if $k=2$, otherwise be equal to an algebraic number times a product of
values of the gamma function at rational arguments by the
Lerch, Chowla--Selberg formula.

The rest of the work consists simply in specializing the above general
argument to specific modular functions $F$ and $h$ and specific CM points
$\tau$.

\smallskip

{\bf Remark.} Under this modular interpretation the existence of these $1/\pi$
formulas is due exclusively to the existence of the modularity-preserving
nonholomorphic modification $D^*$ of the differential operator $D$ seen above,
which involves $1/\pi$.

\subsection{First Special Case: Level $1$}

We first consider modular forms on the full modular group. It is known at
least since Klein--Fricke that we have the hypergeometric representation
$$E_4^{1/4}={}_2F_1(1/12,5/12;1;1/J_1)\;,$$ where
$E_k=1-B_k/(2k)\sum_{n\ge1}\sigma_{k-1}(n)q^n$ and $J_1(\tau)=j(\tau)/1728$.
Thanks to the Clausen identity, we deduce that
$$E_4^{1/2}={}_3F_2(1/2,1/6,5/6;1,1;1/J_1)\;.$$ This is modular of weight $2$,
so we apply the above proposition to $F=E_4^{1/2}$ and $h=1/J_1$.
We compute that $D(h)/h=-D(j)/j=E_6/E_4$, $D(F)/F=(1/6)(E_2-E_6/E_4)$,
hence $D^*(F)/F=(1/6)(E_2^*-E_6/E_4)$ where $E_2^*=E_2-3/(\pi\Im(\tau))$,
so $G^*=-(1/6)(1-E_2^*E_4/E_6)$, so the general identity specializes to
$$\sum_{n\ge0}\left(6n+1-\dfrac{E_2^*E_4}{E_6}(\tau)\right)\dfrac{R_n(2)R_n(6)}{J_1(\tau)^n}=\dfrac{3}{\pi\Im(\tau)}\dfrac{E_4^{3/2}}{E_6}(\tau)\;,$$
an identity due to the Chudnovsky brothers.

For comparison with higher levels, we set $s_1=1/6$, so that for instance
$E_4^{1/4}={}_2F_1(s_1/2,(1-s_1)/2;1;1/J_1)$.

\subsection{Special Cases: Levels $2$ and $3$}

There is no difference for higher levels compared to level $1$, apart from the
need to give explicitly the modular functions used and the hypergeometric
identities.

As it happens, levels $2$ and $3$ can be treated together. For $N=2$ and $3$
set
\begin{align*}F_2(\tau)&=\dfrac{NE_2(N\tau)-E_2(\tau)}{N-1}\;,\quad F_4(\tau)=\dfrac{N^2E_4(N\tau)-E_4(\tau)}{N^2-1}\;,\\
J_N(\tau)&=\dfrac{F_2^4}{F_2^4-F_4^2}\;,\text{\quad and\quad} P_2(\tau)=\dfrac{NE_2(N\tau)+E_2(\tau)}{N+1}\;.\end{align*}
The hypergeometric identity is
$$F_2^{1/2}={}_2F_1(s_N/2,(1-s_N)/2;1;1/J_N)\;,\text{\quad with\quad}s_N=(N+1)/12\;,$$
so by Clausen $$F_2={}_3F_2(1/2,s_N,1-s_N;1,1;1/J_N)\;.$$ We apply the
proposition to $F=F_2$ and $h=1/J_N$. We compute that
$D(h)/h=F_4/F_2$, $D(F)/F=s_N(P_2-F_4/F_2)$, hence
$D^*(F)/F=s_N(P_2^*-F_4/F_2)$ with $P_2^*=P_2-6/((N+1)\pi\Im(\tau))$,
so $G^*=-s_N(1-P_2^*F_2/F_4)$, and since $6/(N+1)=1/(2s_N)$, the general
identity specializes to the two identities
\begin{align*}
\sum_{n\ge0}\left(4n+1-\dfrac{P_2^*F_2}{F_4}(\tau)\right)\dfrac{R_n(2)R_n(4)}{J_2(\tau)^n}&=\dfrac{2}{\pi\Im(\tau)}\dfrac{F_2^2}{F_4}(\tau)\;,\\
\sum_{n\ge0}\left(3n+1-\dfrac{P_2^*F_2}{F_4}(\tau)\right)\dfrac{R_n(2)R_n(3)}{J_3(\tau)^n}&=\dfrac{3}{2\pi\Im(\tau)}\dfrac{F_2^2}{F_4}(\tau)\;.\end{align*}

\subsection{Special Case: Level $4$}

Here we set
\begin{align*}F_2(\tau)&=\dfrac{4E_2(4\tau)-E_2(\tau)}{3}\;,\quad G_2(\tau)=4E_2(4\tau)-4E_2(2\tau)+E_2(\tau)\\
  J_4&=\dfrac{F_2^2}{F_2^2-G_2^2}\;,\text{\quad and\quad}P_2(\tau)=E_2(2\tau)\;.\end{align*}
The hypergeometric identity is again
$$F_2^{1/2}={}_2F_1(s_N/2,(1-s_N)/2;1;1/J_4)\;,\text{\quad with\quad}s_4=1/2\;,$$
so by Clausen $F_2={}_3F_2(1/2,s_4,1-s_4;1,1;1/J_4)$. We apply the
proposition to $F=F_2$ and $h=1/J_4$. We compute that
$D(h)/h=G_2$, $D(F)/F=(P_2-G_2)/3$, hence $D^*(F)/F=(P_2^*-G_2)/3$ with
$P_2^*=P_2-3/(2\pi\Im(\tau))$, so $G^*=-(1-P_2^*/G_2)/3$, hence the
general identity specializes to
$$\sum_{n\ge0}\left(3n+1-\dfrac{P_2^*}{G_2}(\tau)\right)\dfrac{R_n(2)^3}{J_4(\tau)^n}=\dfrac{3}{2\pi\Im(\tau)}\dfrac{F_2}{G_2}(\tau)\;.$$

\section{The Basic List of Rational Hypergeometric $1/\pi$ Formulas}

From the above four specializations it is now immediate to obtain as many
hypergeometric $1/\pi$ formulas as we like. To obtain such formulas which are
\emph{rational} in the above sense, we first need $J_N(\tau)$ to be rational.
This trivially implies that the coefficients of $1/(\pi\Im(\tau))$ on the
right-hand side of the formulas are square roots of rational numbers.
Thus, if we want the coefficient of $1/\pi$ to be algebraic, we need both
$J_N(\tau)$ rational and $\Im(\tau)$ algebraic, and a transcendence theorem 
(well-known for $N=1$, but proved similarly for $N>1$) implies that
$\tau$ is a CM point, i.e., of the form $(a+\sqrt{D})/b$ with
$D<0$ and $a$, $b$ integral. In turn this implies (less trivially)
that the other coefficients involved will be square roots of a rational number.

\medskip

The following table summarizes the results obtained in this way: each
formula is of the form
$\sum_{n\ge0}P(n)H_N(n)/a^n=\sqrt{k}/\pi$,
where the function $H_N(n)=(1/2)_n(s_N)_n(1-s_N)_n/n!^3$ is the coefficient of
$x^n$ in ${}_3F_2(1/2,s_N,1-s_N;1,1;x)$, so that
$H_1(n)=R_n(2)R_n(6)$, $H_2(n)=R_n(2)R_n(4)$, $H_3(n)=R_n(2)R_n(3)$,
and $H_4(n)=R_n(2)^3$. For uniqueness, we always choose $P$ with
content $1$ and positive leading coefficient, $k$ is given in factored
form, and the square root of $k$ is always the positive one. For future
reference, we assign a number from 1 to 36 to each formula.

\vfill\eject

$$\sum_{n\ge0}P(n)\dfrac{H_N(n)}{a^n}=\dfrac{\sqrt{k}}{\pi}$$

\bigskip

\centerline{
\begin{tabular}{|c||c|c|c|c|c|c|}
\hline
\# & $N$ & $H_N$ & $\tau$ & $a=J_N(\tau)$ & $P$ & $k$ \\
\hline\hline
1 & $1$ & $2\cdot6$ & $(1+\sqrt{-7})/2$ & $-2^{-6}5^3$ & $63x+8$ & $3\cdot5^3$ \\
2 & $1$ & $2\cdot6$ & $(1+\sqrt{-11})/2$ & $-2^93^{-3}$ & $154x+15$ & $2^{11}$ \\
3 & $1$ & $2\cdot6$ & $(1+\sqrt{-19})/2$ & $-2^9$ & $342x+25$ & $2^{11}\cdot3$ \\
4 & $1$ & $2\cdot6$ & $(1+\sqrt{-27})/2$ & $-2^93^{-2}5^3$ & $506x+31$ & $2^{11}\cdot3^{-3}\cdot5^3$ \\
5 & $1$ & $2\cdot6$ & $(1+\sqrt{-43})/2$ & $-2^{12}5^3$ & $5418x+263$ & $2^{14}\cdot3^{-1}\cdot5^3$ \\
6 & $1$ & $2\cdot6$ & $(1+\sqrt{-67})/2$ & $-2^95^311^3$ & $261702x+10177$ & $2^{11}\cdot3\cdot5^3\cdot11^3$ \\
7 & $1$ & $2\cdot6$ & $(1+\sqrt{-163})/2$ & $-2^{12}5^323^329^3$ & $545140134x+13591409$ & $2^{14}\cdot3\cdot5^3\cdot23^3\cdot29^3$ \\
8 & $1$ & $2\cdot6$ & $\sqrt{-2}$ & $3^{-3}5^3$ & $28x+3$ & $5^3$ \\
9 & $1$ & $2\cdot6$ & $\sqrt{-3}$ & $2^{-2}5^3$ & $11x+1$ & $2^{-2}\cdot3^{-1}\cdot5^3$ \\
10 & $1$ & $2\cdot6$ & $\sqrt{-4}$ & $2^{-3}11^3$ & $63x+5$ & $2^{-4}\cdot3\cdot11^3$ \\
11 & $1$ & $2\cdot6$ & $\sqrt{-7}$ & $2^{-6}5^317^3$ & $133x+8$ & $2^{-2}\cdot3^{-5}\cdot5^3\cdot17^3$ \\
12 & $2$ & $2\cdot4$ & $(1+\sqrt{-5})/2$ & $-2^2$ & $20x+3$ & $2^6$ \\
13 & $2$ & $2\cdot4$ & $(1+\sqrt{-7})/2$ & $-2^{-8}3^47^2$ & $65x+8$ & $3^4\cdot7$ \\
14 & $2$ & $2\cdot4$ & $(1+\sqrt{-9})/2$ & $-2^43$ & $28x+3$ & $2^8\cdot3^{-1}$ \\
15 & $2$ & $2\cdot4$ & $(1+\sqrt{-13})/2$ & $-2^23^4$ & $260x+23$ & $2^6\cdot3^4$ \\
16 & $2$ & $2\cdot4$ & $(1+\sqrt{-25})/2$ & $-2^63^45$ & $644x+41$ & $2^{10}\cdot3^4\cdot5^{-1}$ \\
17 & $2$ & $2\cdot4$ & $(1+\sqrt{-37})/2$ & $-2^23^47^4$ & $21460x+1123$ & $2^6\cdot3^4\cdot7^4$ \\
18 & $2$ & $2\cdot4$ & $\sqrt{-1}$ & $2^{-5}3^4$ & $7x+1$ & $2^{-2}\cdot3^4$ \\
19 & $2$ & $2\cdot4$ & $\sqrt{-6}/2$ & $3^2$ & $8x+1$ & $2^2\cdot3$ \\
20 & $2$ & $2\cdot4$ & $\sqrt{-10}/2$ & $3^4$ & $10x+1$ & $2^{-3}\cdot3^4$ \\
21 & $2$ & $2\cdot4$ & $\sqrt{-18}/2$ & $7^4$ & $40x+3$ & $3^{-3}\cdot7^4$ \\
22 & $2$ & $2\cdot4$ & $\sqrt{-22}/2$ & $3^411^2$ & $280x+19$ & $2^2\cdot3^4\cdot11$ \\
23 & $2$ & $2\cdot4$ & $\sqrt{-58}/2$ & $3^811^4$ & $26390x+1103$ & $2^{-3}\cdot3^8\cdot11^4$ \\
24 & $3$ & $2\cdot3$ & $(3+\sqrt{-27})/6$ & $-2^43^{-2}$ & $5x+1$ & $2^4\cdot3^{-1}$ \\
25 & $3$ & $2\cdot3$ & $(3+\sqrt{-51})/6$ & $-2^4$ & $51x+7$ & $2^4\cdot3^3$ \\
26 & $3$ & $2\cdot3$ & $(3+\sqrt{-75})/6$ & $-2^4 5$ & $9x+1$ & $2^4\cdot3\cdot5^{-1}$ \\
27 & $3$ & $2\cdot3$ & $(3+\sqrt{-123})/6$ & $-2^{10}$ & $615x+53$ & $2^{10}\cdot3^3$ \\
28 & $3$ & $2\cdot3$ & $(3+\sqrt{-147})/6$ & $-2^4 3^3 7$ & $165x+13$ & $2^4\cdot3^6\cdot7^{-1}$ \\
29 & $3$ & $2\cdot3$ & $(3+\sqrt{-267})/6$ & $-2^4 5^6$ & $14151x+827$ & $2^4\cdot3^3\cdot5^6$ \\
30 & $3$ & $2\cdot3$ & $\sqrt{-6}/3$ & $2$ & $6x+1$ & $3^3$ \\
31 & $3$ & $2\cdot3$ & $\sqrt{-12}/3$ & $2^{-1}3^3$ & $15x+2$ & $2^{-4}\cdot3^6$ \\
32 & $3$ & $2\cdot3$ & $\sqrt{-15}/3$ & $2^{-2}5^3$ & $33x+4$ & $2^{-2}\cdot3^3\cdot5^2$ \\
33 & $4$ & $2^3$ & $(1+\sqrt{-2})/2$ & $-1$ & $4x+1$ & $2^2$ \\
34 & $4$ & $2^3$ & $(1+\sqrt{-4})/2$ & $-2^3$ & $6x+1$ & $2^3$ \\
35 & $4$ & $2^3$ & $\sqrt{-3}/2$ & $2^2$ & $6x+1$ & $2^4$ \\
36 & $4$ & $2^3$ & $\sqrt{-7}/2$ & $2^6$ & $42x+5$ & $2^8$ \\
\hline
\end{tabular}}

\bigskip

\centerline{Rational hypergeometric formulas for $1/\pi$}

\vfill\eject

Note that a generalization of the proof of the class number $1$ problem for
imaginary quadratic fields \emph{proves} that the values for $a$ listed above
(together with the values for the divergent series that we will give below)
are the \emph{only} rational values of $J_N(\tau)$ at CM arguments $\tau$,
outside of the values $0$ and $1$ which cannot be used.

Note that we do not claim that we have found all possible rational $1/\pi$
formulas, but only that, as far as we can tell, no other such formula exists
in the literature, and a rather long search using linear dependence algorithms
did not find any additional ones, except from trivial modifications coming
from the fact that ${}_2F_1$ is solution of a linear differential equation of
order $2$. For instance, we have the following formulas, which are trivially
equivalent to the last three formulas of the above list:
\begin{align*}
\sum_{n\ge0}(6n^3+n^2)\dfrac{R_n(2)^3}{(-8)^n}&=-\dfrac{\sqrt{2}/6}{\pi}\;,\quad\sum_{n\ge0}(2n^3-n^2)\dfrac{R_n(2)^3}{4^n}=\dfrac{1/3}{\pi}\;,\\
\sum_{n\ge0}(210n^3-5n^2+n)\dfrac{R_n(2)^3}{64^n}&=\dfrac{4/3}{\pi}\;.
\end{align*}

Note that the same method allows us to find \emph{divergent} series
because $|a|<1$:

\bigskip

\centerline{
\begin{tabular}{|c||c|c|c|c|c|c|c|}
\hline
\# & $N$ & $H_N$ & $\tau$ & $a=J_N(\tau)$ & $P$ & $k$ & sign\\
\hline\hline
37 & $2$ & $2\cdot4$ & $(-1+\sqrt{-3})/2$ & $-2^{-4}3^2$ & $5x+1$ & $3$ & $+$ \\
38 & $2$ & $2\cdot4$ & $(1+\sqrt{-7})/4$ & $2^{-8}3^4$ & $35x+8$ & $-2^23^4$ & $-$ \\
39 & $3$ & $2\cdot3$ & $(2+\sqrt{-2})/6$ & $2\cdot3^{-3}$ & $10x+3$ & $-2^25^2$ & $-$ \\
40 & $3$ & $2\cdot3$ & $(1+\sqrt{-11})/6$ & $2^43^{-3}$ & $11x+3$ & $-2^43^2$ & $+$ \\
41 & $3$ & $2\cdot3$ & $(3+\sqrt{-15})/6$ & $-2^{-2}$ & $15x+4$ & $3^3$ & $+$ \\
42 & $4$ & $2^3$ & $(1+\sqrt{-1})/2$ & $-2^{-3}$ & $3x+1$ & $1$ & $+$ \\
43 & $4$ & $2^3$ & $(3+\sqrt{-7})/8$ & $2^{-6}$ & $21x+8$ & $-2^4$ & $-$ \\
44 & $4$ & $2^3$ & $(1+\sqrt{-3})/4$ & $2^{-2}$ & $3x+1$ & $-2^2$ & $-$ \\
\hline\hline
\end{tabular}}

\bigskip

The values of $k$ given in this table are those coming from the general
formula, but correspond to the values obtained from the analytic
continuation of the hypergeometric series only when $a<0$. On the other hand,
we will see below that they all lead to supercongruences, as well as so-called
upside-down formulas. Here the sign of the square root of $k$ can vary, so is
indicated in the last column with respect to the principal determination.

\section{Additional Consequences of Proposition \ref{prop:main}}

In the previous section, we have only used Proposition \ref{prop:main}
for some very specific pairs $(F,h)$ which lead to rational hypergeometric
formulas for $1/\pi$. It is evidently possible to use it for other pairs:
in particular we could use other subgroups of $\G$, and in particular the
groups $\G_0^*(N)$ for $N>4$, or still the same subgroups that we have
already considered, but with different $F$ (since the subgroups
$\G_0^*(N)$ for $N\le4$ all have genus $0$, there is not much point in
changing the function $h$, since the resulting formulas could also be obtained
by standard hypergeometric identities).

Using $\G_0^*(N)$ for $N>4$ will not lead to identities involving
hypergeometric functions, but for instance to more general functions called
\emph{Heun functions}.

Using different functions $F$ does give additional formulas.
We simply give two examples, without proof since they are once again
direct applications of Proposition \ref{prop:main}. These are examples in
level $1$, so we keep $h=J_1=j/1728$.

First, we choose $F=E_4^{1/4}$, and we find the general formula
$$\sum_{n\ge0}\left(12n+1-\dfrac{E_2^*E_4}{E_6}(\tau)\right)\dfrac{(1/12)_n(5/12)_n/n!^2}{J_1(\tau)^n}=\dfrac{3}{\pi\Im(\tau)}E_4(\tau)^{-1/4}\dfrac{E_4^{3/2}}{E_6}(\tau)\;.$$
Specializing to $\tau=\sqrt{-3}$ and $\tau=\sqrt{-4}$, and as mentioned
above using the Chowla--Selberg formula to compute the values of
$E_4(\tau)$, we obtain the following identities:
\begin{align*}
\sum_{n\ge0}(22n+1)\dfrac{(1/12)_n(5/12)_n}{n!^2}\dfrac{1}{(125/4)^n}&=\dfrac{(2^{4/3}5^{5/4}/3)\pi}{\G(1/3)^3}=\dfrac{2^{1/3}5^{5/4}3^{-1/2}}{B(1/3,1/3)}\\
\sum_{n\ge0}(126n+5)\dfrac{(1/12)_n(5/12)_n}{n!^2}\dfrac{1}{(11/2)^{3n}}&=\dfrac{2^{1/2}3^{1/4}11^{5/4}\pi^{1/2}}{\G(1/4)^2}=\dfrac{2^{1/2}3^{1/4}11^{5/4}}{B(1/4,1/4)}\;,\end{align*}
where $B(a,b)=\G(a)\G(b)/\G(a+b)$ is the beta function.

Choosing instead $F=E_6^{1/6}$ leads to the following general formula:
$$\sum_{n\ge0}\left(12n+1-\dfrac{E_2^*E_6}{E_4^2}(\tau)\right)\dfrac{(1/12)_n(7/12)_n/n!^2}{(1-J_1(\tau))^n}=\dfrac{3}{\pi\Im(\tau)}E_4(\tau)^{-1/4}\left(\dfrac{E_4^{3/2}}{E_6}(\tau)\right)^{-7/6}\;,$$
and specializing to the same values of $\tau$ gives
\begin{align*}
\sum_{n\ge0}(150n+7)\dfrac{(1/12)_n(7/12)_n}{n!^2}\dfrac{1}{(-121/4)^n}&=\dfrac{2^{4/3}11^{7/6}\pi}{\G(1/3)^3}=\dfrac{2^{1/3}3^{1/2}11^{7/6}}{B(1/3,1/3)}\\
\sum_{n\ge0}(726n+29)\dfrac{(1/12)_n(7/12)_n}{n!^2}\dfrac{1}{(-1323/8)^n}&=\dfrac{2^{1/2}3^{5/2}7^{7/6}\pi^{1/2}}{\G(1/4)^2}=\dfrac{2^{1/2}3^{5/2}7^{7/6}}{B(1/4,1/4)}\;.\end{align*}

\section{Generalization I: Rational Hypergeometric Formulas for $1/\pi^c$}

\emph{Finding} rational hypergeometric identities for $1/\pi^c$ is easily
done using linear dependence algorithms based on the LLL algorithm.
\emph{Proving} them is more difficult: among the methods used is the WZ method,
but it is not the only one. In fact some of the identities (and the three
known identities for $c\ge3$) are still conjectural.
\emph{Explaining} them, as we have done for $1/\pi$ formulas has only
started to be done in a recent paper by Dembel\'e et al. \cite{DPVZ}: the
$1/\pi^2$ formulas are linked to Asai $L$-functions attached to Hilbert modular
forms for real quadratic fields. However, this apparently still does not prove
all of them.

In the following table, we list all known formulas, as before coding $H$ as
$\prod_{d\in I} d^{v_d}$, meaning that
$H_n=\prod_{d\in I}R_n(d)^{v_d}$, and the formula is of the form
$$\sum_{n\ge0}P(n)\dfrac{H_n}{a^n}=\dfrac{\sqrt{k}}{\pi^c}\;.$$

\bigskip

\centerline{
\begin{tabular}{|c||c|c|c|c|c|}
\hline
\# & $c$ & $H$ & $a$ & $P$ & $k$ \\
\hline\hline
1 & $2$ & $2^5$ & $-2^2$ & $20x^2+8x+1$ & $2^6$ \\
2 & $2$ & $2^5$ & $-2^{10}$ & $820x^2+180x+13$ & $2^{14}$ \\
3 & $2$ & $2^3\cdot3$ & $2^63^{-3}$ & $74x^2+27x+3$ & $2^8\cdot3^2$ \\
4 & $2$ & $2^3\cdot4$ & $2^4$ & $120x^2+34x+3$ & $2^{10}$ \\
5 & $2$ & $2\cdot3\cdot4$ & $-2^43$ & $252x^2+63x+5$ & $2^8\cdot3^2$ \\
6 & $2$ & $2\cdot3\cdot6$ & $-2^{12}3^{-6}$ & $1930x^2+549x+45$ & $2^{14}\cdot3^2$ \\
7 & $2$ & $2\cdot3\cdot6$ & $-2^{12}5^3$ & $5418x^2+693x+29$ & $2^{14}\cdot5$ \\
8 & $2$ & $2\cdot3\cdot6$ & $3^{-6}5^6$ & $532x^2+126x+9$ & $2^{-4}\cdot3^2\cdot5^6$ \\
9 & $2$ & $2\cdot4\cdot6$ & $-2^{10}$ & $1640x^2+278x+15$ & $2^{16}\cdot3^{-1}$ \\
10 & $2$ & $2\cdot 8$ & $7^4$ & $1920x^2+304x+15$ & $2^6\cdot7^3$ \\
\hline\hline
11 & $3$ & $2^7$ & $2^6$ & $168x^3+76x^2+14x+1$ & $2^{10}$ \\
\hline\hline
12 & $4$ & $2^5\cdot3\cdot4$ & $-2^83^{-3}$ & $4528x^4+3180x^3+972x^2+147x+9$ & $2^{16}\cdot3^2$ \\
13 & $4$ & $2^74$ & $2^{12}$ & $43680x^4+20632x^3+4340x^2+466x+21$ & $2^{22}$ \\
\hline
\end{tabular}}

\bigskip

\centerline{Rational hypergeometric formulas for $1/\pi^c$}

\bigskip

Most formulas for $1/\pi^2$ were found by the second author, the formula for
$1/\pi^3$ was found by B.~Gourevitch and the two formulas for $1/\pi^4$ were
found by Y.~Zhao and J.~Cullen respectively.

\medskip

We can also find in the literature the following divergent formulas
for $1/\pi^c$:

\bigskip

\centerline{
\begin{tabular}{|c||c|c|c|c|c|c|}
\hline
\# & $c$ & $H$ & $a$ & $P$ & $k$ & sign \\
\hline\hline
14 & $2$ & $2^5$ & $-2^{-2}$ & $10x^2+6x+1$ & $2^4$ & $+$ \\
15 & $2$ & $2^5$ & $-2^{-10}$ & $205x^2+160x+32$ & $2^8$ & $+$ \\
16 & $2$ & $2^33$ & $-3^{-3}$ & $28x^2+18x+3$ & $2^23^2$ & $+$ \\
17 & $2$ & $2\cdot5$ & $-2^85^{-5}$ & $483x^2+245x+30$ & $2^85^2$ & $+$ \\
18 & $2$ & $2\cdot3\cdot4$ & $-2^43^{-3}$ & $172x^2+75x+9$ & $2^83^2$ & $+$ \\
\hline\hline
19 & $3$ & $2^7$ & $2^{-6}$ & $21x^3+22x^2+8x+1$ & $-2^23^2$ & $-$ \\
20 & $3$ & $2^53$ & $2^23^{-3}$ & $92x^3+84x^2+27x+3$ & $-2^83^2$ & $-$ \\
\hline\hline
21 & $4$ & $2^55$ & $-2^{10}5^{-5}$ & $5532x^4+5600x^3+2275x^2+425x+30$ & $2^{16}5^2$ & $+$ \\
\hline\hline
\end{tabular}}

\section{Generalization II: Taylor Expansions of $1/\pi^c$ Formulas}

\subsection{Taylor Expansions of $1/\pi$ Formulas}

Following ideas of the second author \cite{Gui5}, \cite{Gui8},
\cite{Gui9}, we are going to generalize the above
formulas by considering them as constant terms of Taylor expansions. Note that
for any $y$ we can define $(x)_{n+y}=\G(x+n+y)/\G(x)$, hence since $n!=(1)_n$:
$$R_{n+x}(d)=\prod_{\substack{1\le i\le d\\\gcd(i,d)=1}}\dfrac{(i/d)_{n+x}}{(1)_{n+x}}\;.$$
Thus $H_{n+x}$ makes sense, so we could define the generalized sum as
$$\sum_{n\ge0}P(n+x)\dfrac{H_{n+x}}{a^{n+x}}\;.$$
However, when $a<0$ the factor $a^x$ introduces parasitic imaginary terms,
so we prefer to define
$$S(H,a,P;x)=\sum_{n\ge0}P(n+x)\dfrac{H_{n+x}}{\sign(a)^n|a|^{n+x}}=
\sum_{n\ge0}P(n+x)\dfrac{H_{n+x}}{a^n|a|^x}\;,$$
which is equal to $\sign(a)^x$ times the previous one.
Note that this is not the only possible normalization. We could also
shift all the Pochhammer indices by $x$ instead of shifing $n$. In all cases,
this would give the above series multiplied by a quotient of products of gamma
functions involving $x$, so the transformation from one to the other is
immediate.

It is clear that $S(H,a,P;x+1)=\sign(a)(S(H,a,P;x)-P(x)H_x/|a|^x)$, so we may
assume if necessary that $x\in[0,1[$.

In view of the existing literature, we can ask at least two questions: first,
give (at least the initial terms of) the power series expansion of
$S(H,a,P;x)$ around $x=0$. Second, give the value of $S(H,a,P;1/2)$.

One observes that if the value of the sum is $a_0\sqrt{-D}/\pi$ with
$a_0\in\Q^*$ and $D$ a negative fundamental discriminant, the expansion is
always of the form
$$S(H,a,P;x)=a_0|D|\left(\dfrac{\sqrt{-D}}{\pi}+0x-a_2|D|L(D,1)x^2-a_3D^2L(D,2)x^3+O(x^4)\right)\;,$$
with the $a_i$ rational, and where we write $L(D,m)$ for
$\sum_{n\ge1}\lgs{D}{n}/n^m$ (of course, since $D<0$ we have
$a_2|D|L(D,1)=a'_2\sqrt{-D}\pi$ for some rational $a'_2$).
Note that, as mentioned above, it is in principle possible to compute the
coefficient of $x^4$, but by laziness we have done so only for cases (33),
(35), and (36), see below.

The following table uses the same numbering of the formulas as that given
above; the column $C_3$ is related to supercongruences and will be explained
below:

\vfill\eject

\begin{align*}
S(H,a,P;x)&=a_0|D|\left(\dfrac{\sqrt{-D}}{\pi}+0x-a_2|D|L(D,1)x^2-a_3D^2L(D,2)x^3+O(x^4)\right)\;,\\
S_p(H,a,P)&\equiv P(0)\leg{D}{p}p+C_3L(D,3-p)p^3\pmod{p^4}\;.
\end{align*}

\medskip

\centerline{
\begin{tabular}{|c||c|c|c|c|c|c||c|}
\hline
\# & $D$ & $a_0$ & $a_2$ & $a_3$ & $\pi^2S(H,a,P;1/2)/(a_0|D|\sqrt{-D})$ & $C_3$ \\
\hline\hline
1 & $-15$  & $1/3$     & $3/4$  & $1/2$ & $\log(3^3/5)$ & $20$ \\
2 & $-8$   & $2$    & $7/2$  & $4$   & $\log(2)$     & $15$ \\
3 & $-24$  & $2/3$    & $15/4$ & $2$   & $\log(2^5/3^3)$ & $5/2$ \\
4 & $-120$ & $2/27$  & $23/8$ & $1/3$   & $\log(3^3.5/2^7)$ & $5/12$ \\
5 & $-15$  & $128/9$ & $39/4$ & $12$  & $\log(2^23^9/5^7)$ & $5/64$ \\
6 & $-1320$ & $2/3$  & $63/16$ & $1/26$ & $\log(2^{13}11^5/(3^35^{11}))$ & $5/104$ \\
7 & $-40020$ & $16/3$ & $159/128$ & $1/11560$ & $\log(3^{21}5^{13}29^5/(2^{38}23^{11}))$ & $5/36992$ \\
8 & $-20$    & $1/8$  & $1$   & $1/2$ & $2\asin(3/5)$ & $-15/2$ \\
9 & $-15$    & $1/18$ & $2$   & $3/2$ & $2\asin(7/5^2)$ & $-5/8$ \\
10 & $-132$  & $1/96$ & $3/2$ & $1/8$ & $2\asin(41/(3^311))$ & $-5/4$ \\
11 & $-255$  & $1/162$ & $1$  & $1/12$ & $2\asin(4207/(5^417^2))$ & $-5/81$ \\
12 & $-4$    & $1$    & $3$   & $4$ & $2\log(2)$ & $6$ \\
13 & $-7$    & $9/7$  & $5/2$ & $5/2$ & $2\log((88+13\sqrt{7})/3^4)$ & $20/3$ \\
14 & $-3$    & $16/9$ & $21/2$   & $20$  & $(3/2)\log(3^3/2^4)$ & $15/8$ \\
15 & $-4$    & $9$   & $11$  & $20$  & $2\log(3^2/2^3)$ & $10/3$ \\
16 & $-20$   & $36/25$ & $23/4$ & $4$ & $\log(2^{18}/(3^45^5))$ & $1/6$ \\
17 & $-4$    & $441$  & $35$  & $100$  & $2\log(2.3^{10}/7^6)$ & $50/147$ \\
18 & $-4$    & $9/16$  & $2$   & $5/2$ & $2\asin(7/3^2)$ & $-10/3$ \\
19 & $-3$    & $2/3$  & $6$   & $10$ & $\pi/3$ & $-15/8$ \\
20 & $-8$    & $9/64$ & $4$   & $4$  & $2\asin(17/3^4)$ & $-1/3$ \\
21 & $-3$    & $49/27$ & $24$ & $60$ & $2\asin(239/(2.7^4))$ & $-45/392$ \\
22 & $-11$   & $18/11$ & $10$ & $10$ & $2\asin(353/(2/3^8))$ & $-5/24$ \\
23 & $-8$    & $9801/64$ & $28$ & $60$ & $2\asin(8668855388657/(3^811^{12}))$ & $-5/1089$ \\
24 & $-3$    & $4/9$  & $5/2$ & $10/3$ & $\log(3^3/2^2)$ & $5/2$ \\
25 & $-3$    & $4$    & $13/2$ & $10$  & $3\log(4/3)$ & $15/2$ \\
26 & $-15$   & $4/75$ & $7/4$ & $1$    & $\log(3^9/(2^25^5))$ & $1/4$ \\
27 & $-3$    & $32$   & $37/2$ & $40$  & $3\log(2^8/3^5)$ & $15/4$ \\
28 & $-7$    & $108/49$ & $15/2$ & $10$ & $\log(7^7/(2^{10}3^6))$ & $5/18$ \\
29 & $-3$  & $500$  & $85/2$ & $130$ & $3\log(5^6/(2^63^5))$ & $39/50$ \\
30 & $-3$  & $1$  & $2$  & $5/2$ & $(2/3)(\pi-\asin(1633/3^9))$ & $-15/4$ \\
31 & $-4$  & $27/32$ & $4$ & $5$ & $2\asin(329/3^6)$ & $-20/9$ \\
32 & $-3$  & $5/2$ & $8$ & $13$ & $2\asin(239/3^6)$ & $-78/25$ \\
33 & $-4$  & $1/4$ & $1$ & $1$ & $8L(-4,2)/\pi$ & $2$ \\
34 & $-8$ & $1/8$ & $3/2$ & $1$ & $4L(-4,2)/\pi$ & $1$ \\
35 & $-4$ & $1/2$ & $2$ & $2$ & $\pi/2$ & $-2$ \\
36 & $-4$ & $2$  & $6$ & $8$ & $\pi/6$ & $-2$ \\
\hline\hline
\end{tabular}}

\vfill\eject

\bigskip

\subsection{Observations for $1/\pi$}

Concerning the Taylor expansions around $x=0$, note that the
coefficient $a_1$ of $x^1$ always vanishes, and that the numerator
of the first seven $a_2$ are equal to $|D(\tau)|-4$, where $D(\tau)=-7$,
$-11$, $-19$, $-27$, $-43$, $-67$, and $-163$ is the discriminant of the
corresponding $\tau$ (not to be confused with the $D$ occurring in the result).

In addition, we also notice a common pattern for the value at $x=1/2$:

\begin{enumerate}
\item In the value of $S(H,a,P;1/2)$ for $a<0$: with only a few exceptions listed
below, we have $S(H,a,P;1/2)=c_0|D|\sqrt{-D}\log(c_1)/\pi^2$, where
$c_0$ and $c_1$ are rational. The exceptions are as follows:
\begin{itemize}\item In cases (33) and (34) we have
$S(H,a,P;1/2)=c_0|D|\sqrt{-D}L(-4,2)/\pi^3$ with $c_0$ rational.
\item In case (13) $c_1$ is not rational but in the quadratic field
$\Q(\sqrt{7})$.
\end{itemize}
\item In the value of $S(H,a,P;1/2)$ for $a>0$: in all cases we have 
$S(H,a,P;1/2)=c_0|D|\sqrt{-D}\asin(c_1)/\pi^2$ where $c_0$ and
$c_1$ are rational (note that $c|D|\sqrt{-D}/\pi$ is of course of this form,
for instance by choosing $c_1=1$).
\end{enumerate}

It is also possible to guess the coefficient of $x^4$: for instance in
cases (33), (35), and (36), set
$$C_1=\sum_{n\ge1}(-1)^n\dfrac{H_{2n}}{(2n+1)^2}\text{\qquad and\qquad} C_2=\sum_{n\ge1}(-1)^n\dfrac{H_{n}}{(2n+1)^2}\;,$$
where here $H_n=\sum_{1\le j\le n}1/j$ is the $n$th harmonic sum. Then 
\begin{align*}S(H,a,P;x)&=a_0|D|\left(\dfrac{\sqrt{-D}}{\pi}+0x-a_2|D|L(D,1)x^2\right.\\&\left.\phantom{=}-a_3D^2L(D,2)x^3+a_4x^4+O(x^5)\right)\;,\end{align*}
with $D=-4$ and
\begin{align*}
a_4&=(8/3)(50C_1-11C_2-22L(-4,2)\log(2))\;,\\
a_4&=(128/3)(10C_1-C_2-2L(-4,2)\log(2))\;,\\
a_4&=128(22C_1-C_2-2L(-4,2)\log(2))\end{align*}
for cases (33), (35), and (36) respectively.

\subsection{Taylor Expansions of $1/\pi^c$ Formulas for $c\ge2$}

For $c=2$ one observes that if the value of the sum is $a_0\sqrt{D}/\pi^2$ with
$a_0\in\Q^*$ and $D$ a fundamental discriminant, the expansion is always of
the form
\begin{align*}S(H,a,P;x)&=a_0D\left(\dfrac{\sqrt{D}}{\pi^2}+0x-a_2D\sqrt{D}x^2+0x^3\right.\\&\phantom{=}\left.+a_4D^3L(D,2)x^4-a_5a_0D^4L(D,3)x^5+O(x^6)\right)\;,\end{align*}
with the $a_i$ rational (of course, when $D>0$ we have $a_4D^3L(D,2)=a'_4\pi^2D\sqrt{D}$ for some rational $a'_4$).

\begin{align*}
S(H,a,P;x)&=a_0D\left(\dfrac{\sqrt{D}}{\pi^2}+0x-a_2D\sqrt{D}x^2+0x^3\right.\\&\phantom{=}\left.+a_4D^3L(D,2)x^4-a_5a_0D^4L(D,3)x^5+O(x^6)\right)\;,\\
S_p(H,a,P)&\equiv P(0)\leg{D}{p}p^2+C_5L(D,4-p)p^5\pmod{p^6}\;.
\end{align*}

\medskip

\centerline{
\begin{tabular}{|c||c|c|c|c|c|c|c||c|}
\hline
\# & $D$ & $a_0$ & $a_2$ & $a_4$ & $a_5$ & $\pi^3S(1/2)/(a_0D\sqrt{D})$ & $C_5$ \\
\hline\hline
1 & $1$ & $8$ & $1/2$ & $25/4$ & $7$ & $14\z(3)/\pi^2$ & $-7/2$ \\
2 & $1$ & $128$ & $5/2$ & $305/4$ & $7$ & $2\z(3)/\pi^2$ & $-7/2$ \\
3 & $1$ & $48$ & $1/3$ & $4$ & $7/9$ & $2\pi/3$ & $21$ \\
4 & $1$ & $32$ & $1$ & $20$ & $7$ & $\pi/3$ & $21/2$ \\
5 & $1$ & $48$ & $3/2$ & $157/4$ & $91/9$ & $4\log(2^5/3^3)$ & $-91/9$ \\
6 & $1$ & $384$ & $5/6$ & $85/4$ & $7/9$ & $2\log(2)$ & $-315$ \\
7 & $5$ & $128/5$ & $3/2$ & $887/32$ & $21/8$ & $\log(2^{74}5^5/3^{54})$ & $-35/216$ \\
8 & $1$ & $375/4$ & $4/3$ & $40$ & $6944/1125$ & $2\asin(164833/5^8)$ & $1953/50$ \\
9 & $12$ & $32/9$ & $7/24$ & $757/1152$ & $3/16$ & $\log(3^9/2^{14})$ & $-15/2$ \\
10 & $28$ & $1$ & $1/7$ & $31/336$ & $1/21$ & $2\asin(2241857/2^{25})$ & $35/8$ \\
\hline\hline
\end{tabular}}

\bigskip

The last three (one for $1/\pi^3$ and two for $1/\pi^4$) obey completely
similar expansions, but we give them one by one:

\medskip

11:
\begin{align*}S(H,a,P;x)&=32\left(\dfrac{1}{\pi^3}+0x-\dfrac{1}{\pi}x^2+0x^3+(16/3)L(-4,1)x^4+0x^5\right.\\
&\left.\phantom{=}-(8224/45)L(-4,3)x^6+32^2\cdot48L(-4,4)x^7+O(x^8)\right)\;,\\
S(H,a,P;1/2)&=\dfrac{8}{\pi^3}\;,\\
S_p(H,a,P)&\equiv\leg{-4}{p}p^3-6L(-4,5-p)p^7\pmod{p^8}\;.\end{align*}

\medskip

12:
\begin{align*}S(H,a,P;x)&=768\left(\dfrac{1}{\pi^4}+0x-\dfrac{1/2}{\pi^2}x^2+0x^3+(3/8)x^4+0x^5-(147/8)\z(2)x^6\right.\\
&\left.\phantom{=}+0x^7+(471187/1344)\z(4)x^8-3968\z(5)x^9+O(x^{10})\right)\;,\\
S(H,a,P;1/2)&=\dfrac{9216\z(3)}{\pi^7}\;,\\
S_p(H,a,P)&\equiv9p^4-(837/2)\z(6-p)p^9\pmod{p^{10}}\;. 
\end{align*}

\medskip

13:
\begin{align*}S(H,a,P;x)&=2048\left(\dfrac{1}{\pi^4}+0x-\dfrac{2}{\pi^2}x^2+0x^3+(11/3)x^4+0x^5-(908/15)\z(2)x^6\right.\\
&\left.\phantom{=}+0x^7+(53932/7)\z(4)x^8-95232\z(5)x^9+O(x^{10})\right)\;,\\
S(H,a,P;1/2)&=\dfrac{2048/15}{\pi^4}\;,\\
S_p(H,a,P)&\equiv21p^4+(279/4)\z(6-p)p^9\pmod{p^{10}}\;. 
\end{align*}

\bigskip

The observations for $1/\pi^c$ are essentially identical to the case of
$1/\pi$, in particular the coefficients of $x^{2j-1}$ for $1\le j\le c$
vanish.

\section{Generalization III: Upside-Down Series}

For completeness, we list the upside-down series (i.e., with $H_n$ in the
denominator) given in \cite{Gui-Rog}. We do not know if all have been proved,
but probably not all among those with $c>1$.

The general recipe is as follows: if $\sum_{n\ge0}P(n)H_n/a^n=\sqrt{k}/\pi^c$
is a \emph{divergent} or semi-convergent series (i.e., with $|a|<1$ or
$a=-1$) then
$$\sum_{n\ge1}\dfrac{P(-n)}{n^{2c+1}H_n(1/a)^n}=A\cdot L(D,c+1)\;,$$
where $D$ is the fundamental discriminant corresponding to $(-1)^ck$ and
$A\in \Q^*$. Thus, the only new value is that of $A$. However, for the
reader's convenience, we give the list explicitly, the numbering corresponding
to that of the initial series (which is different for $c=1$ and $c\ge2$).

\vfill\eject

$$\sum_{n\ge1}\dfrac{Q(n)}{n^{2c+1}H_nb^n}=A\cdot L(D,c+1)\;.$$

\centerline{
\begin{tabular}{|c||c|c|c|c|c|c|}
\hline
\# & $c$ & $H_N$ & $b=1/a$ & $Q$ & $D$ & $A$\\
\hline\hline
37 & $1$ & $2\cdot4$ & $-2^43^{-2}$ & $5x-1$ & $-3$ & $-45/2$ \\
38 & $1$ & $2\cdot4$ & $2^83^{-4}$ & $35x-8$ & $1$ & $72$ \\
39 & $1$ & $2\cdot3$ & $2^{-1}\cdot3^3$ & $10x-3$ & $1$ & $3$ \\
40 & $1$ & $2\cdot3$ & $2^{-4}3^3$ & $11x-3$ & $1$ & $48$ \\
41 & $1$ & $2\cdot3$ & $-2^2$ & $15x-4$ & $-3$ & $-27$ \\
42 & $1$ & $2^3$ & $-2^3$ & $3x-1$ & $-4$ & $-2$ \\
43 & $1$ & $2^3$ & $2^6$ & $21x-8$ & $1$ & $1$ \\
44 & $1$ & $2^3$ & $2^2$ & $3x-1$ & $1$ & $3$ \\
33 & $1$ & $2^3$ & $-1$ & $4x-1$ & $-4$ & $-16$ \\
\hline\hline
14 & $2$ & $2^5$ & $-2^2$ & $10x^2-6x+1$ & $1$ & $-28$ \\
15 & $2$ & $2^5$ & $-2^{10}$ & $205x^2-160x+32$ & $1$ & $-2$ \\
16 & $2$ & $2^33$ & $-3^3$ & $28x^2-18x+3$ & $1$ & $-14$ \\
17 & $2$ & $2\cdot5$ & $-2^{-8}5^5$ & $483x^2-245x+30$ & $1$ & $-896$ \\
18 & $2$ & $2\cdot3\cdot4$ & $-2^{-4}3^3$ & $172x^2-75x+9$ & $1$ & $-1792$ \\
\hline\hline
19 & $3$ & $2^7$ & $2^6$ & $21x^3-22x^2+8x-1$ & $1$ & $45/4$ \\
20 & $3$ & $2^53$ & $2^{-2}3^3$ & $92x^3-84x^2+27x-3$ & $1$ & $720$ \\
\hline\hline
21 & $4$ & $2^55$ & $-2^{-10}5^5$ & $5532x^4-5600x^3+2275x^2-425x+30$ & $1$ & $-380928$ \\
\hline\hline
\end{tabular}}

\bigskip

\centerline{Upside-down series}

\bigskip

In particular, note that the last formula gives a series for $\z(5)$.

\vfill\eject

We finish by giving an example of a nonrational $1/\pi$ formula, but there
exist almost a hundred involving only quadratic irrationals, listed in
\cite{Ald}:

$$\sum_{n\ge0}\dfrac{R_n(2)^3}{(2+\sqrt{3})^{4n}}(12n+(3-\sqrt{3}))
=\dfrac{(2+\sqrt{3})(4/3)^{1/4}}{\pi}\;.$$

\section{Generalization IV: Supercongruences}

It has been noted long ago by several authors that to \emph{every} $1/\pi^c$
formula (including divergent ones) corresponds a congruence modulo $p$ to a
higher power than could be expected, what is now called a
\emph{supercongruence}.

The main observation is as follows: if there exists a $1/\pi^c$-formula
of the form
$\sum_{n\ge0}P(n)H_n/a^n=\sqrt{k}/\pi^c$, then for all primes $p$ such
that $v_p(a)=v_p(k)=0$ and not dividing any $d$ occuring in $H$ we should
have the following precise supercongruence:
$$\sum_{n=0}^{p-1}P(n)\dfrac{H_n}{a^n}\equiv P(0)\leg{(-1)^c4k}{p}p^c\pmod{p^{2c+1}}\;.$$

For instance, we have the following supercongruences:

\begin{align*}
\sum_{n=0}^{p-1}(154n+15)\dfrac{R_2(n)R_6(n)}{(-8/3)^{3n}}&\equiv 15\leg{-8}{p}p\pmod{p^3}\\
\sum_{n=0}^{p-1}(5418n^2+693n+29)\dfrac{R_2(n)R_3(n)R_6(n)}{(-80)^{3n}}&\equiv 29\leg{20}{p}p^2\pmod{p^5}\end{align*}

The same phenomenon is valid for the \emph{divergent} series for
$1/\pi^c$ that we have given. For instance we have

$$\sum_{n=0}^{p-1}(35n+8)\dfrac{R_2(n)R_4(n)}{(3/4)^{4n}}\equiv 8p\pmod{p^3}$$

However it has been noticed by several authors that these supercongruences
can be refined to a higher power of $p$ (more precisely to a congruence
modulo $p^{2c+2}$ instead of $p^{2c+1}$): the recipe, made precise by
the second author, is simply to replace the $L(D,c+1)$ occurring in the
coefficient of $x^{2c+1}$ of the Taylor expansions by $L(D,c+2-p)$ times a
suitable rational number. In other words,
\begin{align*}S_p(H,a,P)&:=\sum_{n=0}^{p-1}P(n)\dfrac{H_n}{a^n}\\
&\equiv P(0)\leg{(-1)^c4k}{p}p^c+C_{2c+1}p^{2c+1}L(D,c+2-p)\pmod{p^{2c+2}}\;.\end{align*}
The coefficients $C_{2c+1}$ have been given for all the convergent series
in the above tables. The remaining coefficients for the divergent series
are as follows:

For the divergent $1/\pi$ formulas,
$$C_3=(15/4,0,0,0,12,2,0,0)$$
for formulas (37) to (44).

For the divergent $1/\pi^c$ formulas for $c\ge2$,
$$C_{2c+1}=(-7/2,-64,-21/2,-210,-63,0,0,-1395)$$
for formulas (14) to (21).

The coefficients $0$ of course mean that the congruence is valid modulo
$p^{2c+2}$ with no correction terms.

\medskip

We can observe that $C_3$ is almost always divisible by $5$, $C_5$ is almost
always divisible by $7$, and $C_9$ is always divisible by $279=3^2\cdot31$. We
have no explanation for this phenomenon.

\medskip

Henri Cohen, Universit\'e de Bordeaux, LFANT, IMB, U.M.R. 5251 du C.N.R.S,
351 Cours de la Lib\'eration, 33405 Talence Cedex, FRANCE.

\smallskip

Jes\'us Guillera, Universidad de Zaragoza, Departamento de Matem\'aticas,
50009 Zaragoza, SPAIN.

\end{document}